\documentclass{article}

\usepackage{amsmath}
\usepackage{amsthm}
\usepackage{amsfonts}
\title{Holomorphic Diffeomorhisms of Semisimple Homogenous Spaces}
\author{\'Arp\'ad T\'oth
\and Dror Varolin
\thanks{2000 Mathematics Subject Classification 32M05, 32M17, 32M25, 22E46}}
\date{}

\newcommand{\noi}{\noindent}
\newcommand{\pf}{{\noi \it Proof: }}
\newcommand{\rmk}{\noi {\sc Remark: }}

\newcommand{\bs}{\bigskip}

\newcommand{\cc}{{\mathcal C}}

\newcommand{\cl}{{\mathcal L}}

\newcommand{\co}{{\mathcal O}}

\newcommand{\cu}{{\mathcal U}}

\newcommand{\cx}{{\mathcal X}}

\newcommand{\fb}{{\mathfrak b}}

\newcommand{\fg}{{\mathfrak g}}
\newcommand{\fh}{{\mathfrak h}}

\newcommand{\fk}{{\mathfrak k}}
\newcommand{\fl}{{\mathfrak l}}

\newcommand{\fn}{{\mathfrak n}}
\newcommand{\fo}{{\mathfrak o}}
\newcommand{\fp}{{\mathfrak p}}

\newcommand{\fs}{{\mathfrak s}}

\newcommand{\tr}{{\rm{tr}}}

\newcommand{\vp}{\varphi}

\newcommand{\R}{{\mathbb R}}
\newcommand{\p}{{\mathbb P}}

\newcommand{\C}{{\mathbb C}}

\newcommand{\cn}{{\mathbb C} ^n}

\newcommand{\holvec}{\cx _{\co}}
\newcommand{\faff}{\co _{Aff}}
\newcommand{\vaff}{\cx _{Aff}}

\newcommand{\di}{\partial}

\newcommand{\slc}{\fs \fl (2,\C)}

\newcommand{\diff}{{\rm Diff} _{\co}}
\newcommand{\ii}{\sqrt{-1}}

\newcommand{\emb}{\hookrightarrow}
\newcommand{\tensor}{\otimes}

\begin{document}
\maketitle

\newtheorem{thm}{\sc Theorem}[section]
\newtheorem{lem}[thm]{\sc Lemma}
\newtheorem{prop}[thm]{\sc Proposition}
\newtheorem{cor}[thm]{\sc Corollary}
\newtheorem{conj}[thm]{\sc Conjecture}
\newtheorem{defn}[thm]{\sc Definition}
\newtheorem{qn}[thm]{\sc Question}

\section{Introduction}

Complex analytic manifolds whose group of holomorphic diffeomorphisms acts
transitively have been at the center of much research since Klein's
Erlangen program.  Almost as early on, studying these symmetries proved
itself a useful tool in the function theory of several complex
variables. When the homogeneous space can be represented as a bounded
symmetric domain in $\C^n$, the group of holomorphic diffeomorphisms
preserves a natural metric, and is therefore finite dimensional.
Poincar\' e computed these groups in the case of the ball and the
bidisk, and was lead to the failure of the Riemann Mapping Theorem in
higher dimensions.

In this paper we study the the analytic geometry of Stein
manifolds that support a transitive action, by holomorphic
diffeomorphisms, of a complex semisimple Lie group.  Unlike the
aforementioned bounded symmetric domains, these manifolds have
holomorphic diffeomorphism groups that are infinite dimensional, and
thus a different approach is required for their study.
This work is based on certain connections between representation
theory and properties of the entire group of holomorphic
diffeomorphisms that we establish through further development of
our earlier work \cite{tv}.

Past work in the study of affine homogeneous spaces has
concentrated mainly on algebraic automorphisms, invariant theory,
and applications to algebraic geometry.  Even in the simplest case
of $\cn$, the group of holomorphic diffeomorphisms had not been
well understood until the ground-breaking work of Anders\' en and
Lempert \cite{a,al}.  Recall that a shear is a transformation of
$\cn$ that, after an affine change of coordinates, has the form
$(w,z) \mapsto (z, e^{a(z)}w+b(z))$, where $(z,w) \in \C ^{n-1}
\oplus \C$ and $a, b \in \co (\C ^{n-1})$.  The Anders\'
en-Lempert Theorem states that when $n \ge 2$ every holomorphic
diffeomorphism of $\cn$ can be approximated uniformly on compact
sets by compositions of shears.

In a manner similar to $\C ^n$, semisimple homogeneous spaces such
as the smooth affine quadrics $Q_n = SO(n+1,\C)/SO(n,\C) = \left
\{ x_0^2 + ... + x_n ^2 = 1 \right \} \subset \C ^{n+1}$ are very
symmetric and large in the complex analytic sense, and so it is
difficult to say things about the group of holomorphic
diffeomorphisms.

The density property is a notion capturing the intuitive idea that
the group of holomorphic diffeomorphisms is as large as possible.
Let $X$ be a complex manifold, and $\holvec (X)$ the Lie algebra
of holomorphic vector fields on $X$.  We say that $X$ has the {\it
density property} if the Lie subalgebra generated by the complete
vector fields on $X$ is dense in $\holvec (X)$:
$$
\overline {\left < \xi \in \holvec (X)\ ;\ \xi \ {\rm is\
complete} \right >} = \holvec (X).
$$
(See section \ref{background-section} for the definition of
complete vector fields.)

In fact, we establish the density property for a large collection
of semisimple homogeneous spaces. The main result of this paper is
the following theorem.

\begin{thm}\label{adjoint}
If $G$ is a complex semisimple Lie group of adjoint type and $K$ a
reductive subgroup, then the homogeneous space $X=G/K$ has the
density property.
\end{thm}

\noi Recall that a semisimple Lie group is of {\it adjoint} type if
its center is trivial or, equivalently, if the adjoint representation
is faithful.
We were not able to establish the density property for
all semisimple homogeneous spaces; the proof of our main theorem uses
the adjoint condition in a crucial way.

\medskip

The density property is of fundamental significance, and governs
the geometry of the underlying space in various respects, some of
which were discussed in \cite{v1, v2}. By results of \cite{v2}, we
have the following corollaries of Theorem \ref{adjoint}.

\begin{cor}\label{cors}
Let $X$ be as in Theorem \ref{adjoint}, and let $n = {\rm
dim}_{\C}X$.
\begin{enumerate}
\item There is an open cover of $X$ by open subsets each of which is
biholomorphic to $\C^n$, i.e., is a  Fatou-Bieberbach domain.
\item The space $X$ is biholomorphic to one of its proper open
subsets. \item Let $Y$ be a complex manifold with ${\rm dim}(Y) <
{\rm dim}(X)$, such that  and there exists a proper holomorphic
embedding $j: Y \emb X$. Then there exists another proper
holomorphic embedding $j': Y \emb X$ such that for any $\vp \in
{\rm Diff} _{\co} (X)$, $\vp \circ j (Y) \neq j' (Y)$.
\end{enumerate}
\end{cor}

Here and below, ${\rm Diff} _{\co} (X)$ refers to the group of
holomorphic diffeomorphisms of the complex manifold $X$. We avoid
the more standard notation $Aut(X)$, and use the words {\it
holomorphic diffeomorphisms} in place of the more common {\it
automorphisms} for fear that when Lie groups are in the picture,
use of the word {\it automorphism} may cause confusion.

\medskip

Knowing the density property for a given space has applications
not only to the geometry of that space, but also to some
associated spaces on which the density property is not known to
hold.  There are many results of this kind, but for the purpose of
illustration we restrict ourselves to the following theorem, which
is our second main result of the paper.

\begin{thm}\label{lift}
Let $X$ be a semisimple homogeneous space of dimension $n$ and
$Y$ a complex manifold of dimension $< n$
such that there exists a proper holomorphic embedding $j:
Y \emb X$. Then there exists another proper holomorphic embedding
$j': Y \emb X$ such that for any $\vp \in {\rm Diff} _{\co} (X)$,
$\vp \circ j (Y) \neq j' (Y)$.
\end{thm}

It is not known if the spaces $X$ to which Theorem \ref{lift} applies
have the density property.  They are all finite covers of the
spaces in Theorem \ref{adjoint}, and the proof of Theorem
\ref{lift} uses only this fact.  It is an open problem whether
finite covers of Stein manifolds with the density property have
the density property, and vice versa.  If this turns out to be
true, it would greatly simplify the proof of the main result of
\cite{tv}.  In fact, the general relationship between the density
property and covers (finite or not) is still not well understood,
with the most famous open problem being the question of whether or
not there is a Fatou-Bieberbach domain in $\C ^* \times \C ^*$.

\medskip

One of the two central ingredients in the proof of Theorem
\ref{adjoint} is the general notion of {\it shears}, introduced in
\cite{v3}.  This notion, whose key ideas are recalled in section
\ref{background-section}, turns the verification of the density
property on homogeneous spaces into a problem in representation
theory of semisimple Lie groups-- the second central ingredient.

\medskip

The density property was introduced so as to allow generalization
to other complex manifolds of the previously mentioned results of
Anders\'en and Lempert \cite{a,al}.  In the process, this
definition also simplified the proofs of those results, even in
the original setting.  The theorems of Anders\' en and Lempert
began as answers to questions raised by J.-P. Rosay and W. Rudin
in their foundational paper \cite{rr}.  It was also realized (in
part in the paper \cite{al}, and more fully in the joint work of
F. Forstneri\v c and Rosay \cite{fr}) that the density property
could be used in many interesting analytic geometric constructions
in $\cn$. In this regard, there are definitions of the density
property for more general Lie algebras of vector fields (see
\cite{v1}), taking into account whatever additional geometry one
cares about.  For more details on this side of the story, we
recommend the survey papers \cite{f} and \cite{r}.

As it turns out, the notion of the density property transcends
generalization for its own sake, bringing to light new features of
the complex manifolds possessing this property; witness Corollary
\ref{cors} and further results in \cite{v2}.

The density property on a Stein manifold $X$ is a meaningful
statement about the size of the group $\diff (X)$.  Indeed, one
way of measuring the size of this group is to look at its formal
tangent space: the closure of the Lie algebra generated by the
complete vector fields on $X$.  The density property is the
statement that this formal tangent space is as large as possible.

Some additional remarks are of interest.

\noindent {\bf 1.}  Previously, there was some hope that the
"largest" spaces were determined by their topology and the
identical vanishing of the Kobayashi infinitesimal pseudometric.
This was formulated precisely as the following question of
Diederich and Sibony \cite{ds} for the case of $\cn$, $n \ge 2$:

\noi {\it If $X$ is a Stein manifold that is diffeomorphic to $\R
^{2n}$ and has identically vanishing Kobayashi infinitesimal
pseudometric, is $X$ biholomorphic to $\cn$?}

\noi A recent counterexample of Forn\ae ss \cite{fv} answered this
question in the negative.  In view of corollary \ref{cors}.1, we
have the following

\begin{conj}\label{poincare}
If $X$ is a Stein manifold that is diffeomorphic to $\cn$ and has
the density property, then $X$ is biholomorphic to $\cn$.
\end{conj}

\noindent {\bf 2.} From the analytic geometry point of view, Stein
manifolds with the density property seem to share many geometric
properties of rational algebraic varieties, and thus it is natural
to ask the following question.

\begin{qn}\label{moishe}
Does every Stein manifold with the density property, have a
Moishezon compactification?  Is there a compactification that is
bimeromorphic to $\p _n$?
\end{qn}

By combining Corollary \ref{cors}.1 with a famous theorem of
Kodaira \cite{k}, one has that a surface with the density property
that also admits a compactification must be rational.  This is
generally false in higher dimensions.

\medskip

\noi The paper is organized as follows.  In Section
\ref{background-section} we recall various known results that are
needed in the proofs of our main theorems.  Some of these results
are standard, and others were developed in \cite{v3} and
\cite{tv}.  In Section \ref{reduce} we extend results in \cite{tv}
to reduce the study of the density property on complex semisimple
homogeneous spaces of adjoint type to a problem in representation
theory.   Although much of Section \ref{adjoint-proof-section} is
a collection of facts, it is technically the most demanding part
of the paper, and so in Section \ref{examples-section} we show
how, in some specific cases, the density property can be
established without recourse to the technical arguments and
full-blown root system machinery of Section
\ref{adjoint-proof-section}. The homogeneous spaces we use as
examples are ubiquitous, appearing naturally in many mathematical
constructions.  The conclusion of the proof of Theorem
\ref{adjoint} occupies Section \ref{adjoint-proof-section}, using
a case-by-case analysis of the representation theoretic problem
established in Section \ref{reduce}.    The proof of Theorem
\ref{lift} occupies the sixth and final section.

\section{Background Material}\label{background-section}

In this section, we state various results that are needed later.
This is also an opportunity to establish some notation. Since all
of the material in this section is contained elsewhere, we omit
almost all proofs.

\subsection*{Holomorphic vector fields.}

A holomorphic vector field $\xi$ on a complex manifold $X$ is a
holomorphic section of $T^{1,0}_X$, the holomorphic part of the
complexified tangent bundle.  We let $\holvec (X)$ denote the set
of holomorphic vector fields on $X$.  Since $T^{1,0}_X$ is
naturally isomorphic to the real tangent bundle $T_X$, we can
identify $\xi$ with a real vector field  that we continue to
denote by $\xi$.  As such, there is a flow $\vp _{\xi}$ associated
to $\xi$, which is defined on an open subset $\cu$ of $\R \times
X$ containing $\{ 0\} \times X$ in the following way: for $(t,p)
\in \cu$, $\vp _{\xi} ^t (p)  = c(t)$, where $c: (-a(p),b(p)) \to
X$ is the unique maximal solution of the initial value problem
$$\frac{dc}{ dt} = \xi \circ c, \qquad c(0)=p.\eqno{(*)}$$
It follows from general ODE theory that the map $p \mapsto \vp
_{\xi} ^t (p)$ is holomorphic.

We say that $\xi$ is {\it complete} if $\cu = \R \times X$, i.e.,
if for each $p \in X$ one can solve $(*)$ for all $t \in \R$.  In
this case $\{ \vp _{\xi} ^t \ |\ t \in \R \}$ is a one parameter
group of holomorphic diffeomorphisms of $X$.  We say that $\xi$ is
{\it $\C$-complete} if both $\xi$ and $i\xi$ are complete.  Define
the $\C$-flow of $\xi$ to be
$$g_{\xi} ^{s+it}:= \vp _{\xi} ^s \circ \vp _{i\xi} ^t.$$
If $\xi$ is $\C$-complete, then $\{ g _{\xi} ^{\zeta} \ |\ \zeta
\in \C \}$ defines a holomorphic $\C$-action.  In this paper, all
complete vector fields are $\C$-complete, so we shall often drop
the prefix $\C$, and still refer to $g_{\xi}$ as the flow of
$\xi$, even though it is defined for ``complex time''.

With the operation $[\xi ,\eta ] = \xi \eta - \eta \xi $, $\holvec
(X)$ forms a Lie algebra.  We can generate a Lie subalgebra of
$\holvec (X)$ using complete vector fields on $X$. We shall call
any vector field in the closure of this subalgebra {\it completely
generated.}  In general, this subalgebra will not consist of
complete vector fields.  However, completely generated vector
fields have the extraordinary property that their flows can be
approximated (in the locally uniform, and hence $\cc ^k$ topology)
by holomorphic diffeomorphisms of $X$ \cite{v1}.  This result can
be proved as a combination of an approximation method, often
attributed to Euler, and the following formulas. Let $\xi$ and
$\eta$ be vector fields with flows $\vp _{\xi}$ and $\vp _{\eta}$
respectively.
$$\left . \frac{d}{dt} \right |_{t=0} \vp _{\xi}^t \circ \vp
_{\eta} ^t = \xi + \eta $$
$$\left . \frac{d}{dt} \right |_{t=0+} \vp _{\eta} ^{-\sqrt {t}}
\circ \vp _{\xi} ^{-\sqrt {t}} \circ \vp _{\eta} ^{\sqrt {t}}
\circ \vp _{\xi} ^{\sqrt {t}} = [\xi,\eta].$$  The approximation
method of Euler can be stated as follows.
\begin{thm}
Let $X$ be a manifold, and let $\xi$ be a vector field with flow
$\vp _{\xi}$ on $X$.  Suppose $F_t : X \to X$ is a $\cc ^1$ family
of self maps such that $F_0 = id _X$ and
$$\left . \frac{d}{dt} \right | _{t=0} F_t = \xi.$$
Then, in the compact-open topology, one has
$$\lim _{N \to \infty} \underbrace{F_{t/N} \circ ... \circ
F_{t/N}}_{N\ {\rm times}} = \vp _{\xi} ^t.$$
\end{thm}
\noi Putting all of this together, we have the following theorem.
\begin{thm}\label{euler}
Let $\xi$ be a vector field that lies in the closure of the Lie
algebra generated by the complete vector fields on $X$.  Let $K
\subset X$ be compact, and let $t$ be such that the flow of $\xi$
is defined up to time $t$ on $K$.  Then $\vp _{\xi} ^t |K$ can be
approximated, uniformly on $K$, by holomorphic diffeomorphisms of
$X$.
\end{thm}

\subsection*{General shears}

In \cite{v3} the following fundamental proposition was proved.

\begin{prop}\label{dror}
Let $\xi \in \holvec (X)$ be $\C$-complete, and let $f \in \co
(X)$. Then $f \cdot \xi$ is $\C$-complete if and only if $\xi ^2 f
= 0$.
\end{prop}

\noi One is thus naturally lead to study the function spaces
$$I^1(\xi) = \{ f \in \co (X) \ | \ \xi f = 0 \} \qquad {\rm and} \qquad
I^2(\xi ) = \{ f \in \co (X) \ | \ \xi ^2 f = 0 \}$$ consisting of
holomorphic {\it first} and {\it second} integrals respectively.

\begin{defn}\label{shear}
Let $\xi \in \holvec (X)$ be $\C$-complete.  We call $f \cdot \xi$
a $\xi$-shear (resp. $\xi$-overshear ) if $f \in I^1(\xi )$ (resp.
$f \in I^2 (\xi )$).
\end{defn}

\noi We will often refer to $\xi$-shears simply as shears, and
similarly with overshears.

\medskip

\noi \rmk  Let $X=\cn$, let $\xi$ be a constant vector field and
let $f\in I^j(\xi)$ for $j=1,2$.  In this case, it is the time-1
map of $f\xi$, rather than $f\xi$ itself, that is called a shear
and overshear in the literature.  In a sense, our shears are {\it
infinitesimal shears}.

\medskip

For a given complete vector field, the existence of first
integrals is a classical problem and a highly nontrivial matter.
The orbits of the vector field must sit fairly nicely together,
generally speaking.  For the existence of second integrals that
are not first integrals, almost all the orbits must be
biholomorphic to $\C$, and the orbit space must be extremely
regular \cite{v3}. Nevertheless, in the case of complex
(semisimple) Lie groups, many left invariant vector fields have a
lot of first and second integrals.

\subsection*{Semisimple Lie algebras.}

The following notation will be used throughout the paper. We
denote by $\slc$ the three dimensional complex Lie algebra with
basis $\{ E,F,H\}$ satisfying the commutation relations
\begin{equation}\label{sl-2-base}
[H,E]=2E \qquad [H,F]=-2F \qquad [E,F]=H.
\end{equation}

By a representation $\rho$ of  a Lie algebra $\fg$, we mean a
vector space $V$ and a linear map $\rho: \fg \rightarrow End(V)$
satisfying
$$\rho([\xi ,\eta ])=\rho(\xi )\rho( \eta )-\rho(\eta )\rho(\xi ).$$
We will refer to such $V$ as $\fg$-spaces, or $\fg$-modules. Let
$ad: \fg \rightarrow End(\fg)$ be the adjoint representation,
i.e., the representation of $\fg$ on itself given by $ad (\xi)
\eta :=[\xi ,\eta ]$ for all $\xi ,\eta  \in \fg$.

For a general semisimple Lie algebra $\fg$, we fix a Cartan
subalgebra $\fh$, i.e., a commutative subalgebra all of whose
elements can be simultaneously diagonalized in the adjoint
representation, say, and whose dimension is maximal.

Given any representation $V$ of $\fg$, a {\it weight} of $V$ is a
linear functional $\lambda : \fh \to \C$ such that the subspace
$$V_{\lambda} := \left \{ v \in V\ ;\ H v = \lambda(H) v\ {\rm
for\ all}\ H\in \fh \right \}$$ is nonzero.

The roots of $\fh$ are the weights $\alpha$ of the adjoint
representation, and their union forms a {\it root system} $\Phi$.
(See \cite{b} for more on root systems.)  We choose an ordering of
the roots as follows. Let $H \in \fh$ be such that no root
vanishes on $H$. If a root $\alpha$ satisfies $\alpha (H) >0$, we
say that $\alpha$ is positive, and write $\alpha \in \Phi ^+$.
Among the positive roots, there is a subset $\Delta $ of so called
{\it simple roots}, with the property that any positive root can
be written as a linear combination of the simple roots, with
nonnegative integer coefficients.

The simple roots form a basis for $\fh ^*$, the dual vector space
of $\fh$. Moreover, given a root $\alpha$, $-\alpha$ is also a
root, and there are vectors $E_\alpha \in \fg_\alpha, F_\alpha \in
\fg_{-\alpha}$ and $H_\alpha\in \fh$, so that the following hold.

\begin{enumerate}
\item  The Lie algebra spanned by $E_\alpha , F_\alpha$ and
$H_\alpha$ is isomorphic to $\slc$, i.e.,  $E_\alpha , F_\alpha$
and $H_\alpha$ satisfy the commutation relations
(\ref{sl-2-base}).  We define the Lie algebra isomorphism
$$\phi _{\alpha} :\slc \to {\rm span} \left \{E_{\alpha},F_{\alpha},
H_{\alpha} \right \} \quad {\rm by}\quad \phi _{\alpha} \left \{
\begin{array}{c}
E\\F\\H
\end{array}
\right \} = \left \{
\begin{array}{c}
E_{\alpha}\\F_{\alpha}\\H_{\alpha}
\end{array}
\right \}.
$$

\item For all linearly independent roots $\alpha, \beta$,
$$\alpha(H_\beta) = \frac{2 B(\alpha,\beta)}{B(\beta,\beta)},$$
where $B(\alpha,\beta)=\tr ( ad H_\alpha ad H_\beta)$ is the form
on the root system associated to $\fh$ induced by the Killing form
of $\fg$.
\end{enumerate}

\noi The collection of $H_{\alpha}$, $\alpha \in \Phi$, is thus
also a root system $\Phi ^* $, called the dual of $\Phi$.

\subsection*{Homogeneous spaces}

A homogeneous space $X$ is a manifold that has a transitive action by
a Lie group $G$.  In our main theorems, we assume that $G$ is
semisimple and of adjoint type.  Thus we are dealing with
groups having trivial center. \footnote{On occasion, when it is more
convenient to do so,  we might represent
homogeneous spaces in the form $\tilde{G}/K$, where $\tilde{G}$ has
finite center contained in $K$.}  It is known that every
adjoint group is algebraic, i.e. it is isomorphic to a closed subgroup
of $SL_N(\C)$ defined by the vanishing of certain polynomials.  By
Cartan's Theorem $A$, the algebra $\C[G]$, of polynomial functions on
$G$ (where $G$ is viewed as a closed submanifold on End$(\C ^N)$) , is
dense in $\co(G)$. If $K$ is a subgroup of $G$, the quotient manifold
$G/K=\{gK:g\in G\}$ admits a transitive action by $G$, $h: gK
\mapsto hgK$, and every manifold that admits a transitive action
by $G$ arises in this way. By Matsushima's theorem \cite{m}, such
a quotient $G/K$ is Stein if and only if $K$ is reductive.

We now turn to the realization of $G/K$ in the case when $K$ is a
reductive subgroup of $G$.  This is the ``if'' part of the theorem
mentioned above, and is a classical theorem going back to Weyl and
Hilbert.  We use the word {\it reductive} in an analytic sense:
$K$ is reductive if it contains a compact subgroup $K_0$ so that
$\fk = \fk_0+i \fk_0$, where $\fk$ and $\fk_0$ are the Lie
algebras of $K$ and $K_0$ respectively.  With this definition, we
shall need is following result of H. Weyl.  For the proof, see
\cite{mum}.

\begin{thm}[Weyl]\label{realization}
Let $G$ be a semisimple group, and $K$ a reductive subgroup. Then
there is a representation of $G$ in a finite dimensional vector
space $V$, and an element $v \in V$, so that the orbit $X=G\cdot
v$ of $v$ by $G$ is closed, and biholomorphic to $G/K$.
\end{thm}

The realization of $G/K$, via theorem \ref{realization}, as an
orbit in a finite dimensional representation is used in the proof
of our theorems. To see how, let $\rho$ be a linear representation
of $G$ in a vector space $V$.  To every $\xi \in \fg$, we can
associate a complete holomorphic vector field $\vec \xi$ on $V$
defined by
\begin{eqnarray*}
\vec \xi f (p) &=& \left . {\frac{d}{dt}} \right | _{t=0} f (\rho
( e ^{(t\xi )})p) \qquad f \in \co (V),\ p \in V.
\end{eqnarray*}

\noi Evidently $\vec \xi$ is $\C$-complete. If we now let $f$ be a
(holomorphic) defining function for a $G$-invariant closed
submanifold $X$ of $V$, in a neighborhood $U$ (in $V$) of a point
$p \in X$,  then
$$(df \vec \xi )_p  = \vec \xi f (p) = \left . \frac{d}{dt} \right |
_{t=0} f(\rho ( e ^{(t\xi )})p) = 0,$$ since, for $t$ small
enough, $\rho ( e ^{(t\xi )})p \in U \cap X$.  This proves the
following proposition.

\begin{prop}
For every $\xi \in \fg$, the vector field $\vec \xi$ is tangent to
any smooth $G$-orbit in $V$.
\end{prop}

\rmk  {\it Note that Weyl's construction
provides and embedding of $G/K$ into a linear space on which $G$, an
adjoint group, acts. It follows that the weights admitted by this
representation are in the root lattice.}

\section{Reduction to representation theory}\label{reduce}

The first part of the section reviews some of our work in
\cite{tv}. We present a number of technical refinements that will be
used in the proof of Theorem~\ref{adjoint}.

\subsection*{Polynomial vector fields}
We begin by introducing an algebraic version of Definition
\ref{shear}. To this end, let $V$ be $\fg$-space.  The vector
space $V \otimes \fg$ has a natural $\fg$-module structure, given
by
$$(\rho \otimes ad) (\xi )v\tensor \eta =\rho(\xi )v \tensor \eta
+ v \tensor [\xi ,\eta ].$$

\begin{defn} We call $f \otimes \xi$ a $\fg$-shear (resp.
$\fg$-overshear) if $\rho(\xi )f=0$, (resp. $\rho(\xi )^2f = 0$).
\end{defn}

Our goal, then, is to study whether or not, for a given
$\fg$-space $V$, overshears generate $V \tensor \fg$.  The
simplest case to consider is the case $\fg = \slc$, and we did so
in \cite{tv}.  We now recall several results from there, beginning
with the following theorem.

\begin{thm}\label{even}
If $\fg=\slc$ and $V$ is either $\C^2$ or an irreducible module
whose dimension is odd, then $V$ is generated by its
$\slc$-overshears.
\end{thm}

Now, given a $\fg$-space $V$, we can also form the $\fg$-space $\C
[V] \tensor \fg$, where $\C [V]$ is the algebra of polynomial
functions on $V$.  Remarkably, using Theorem \ref{even} one has
the following result.

\begin{thm}\label{algebra} Let $\fg$ be semisimple and
let $V$ be a $\fg$-space. Assume that $V^* \otimes \fg$ is
generated by its $\fg$-overshears. Then $\C[V] \otimes \fg $ is
generated by its $\fg$-overshears.
\end{thm}

\pf We chose a basis $x_i$ of $V^*$ consisting of weight vectors. This
choice defines monomial basis (consisting of weight vectors) of the
finite dimensional vector spaces of polynomials of a given degree. It
suffices to check that for every monomial $f \in \faff (M)$ and $X \in
\fg$ of the form $F_\alpha, E_\alpha,$ or $H_\alpha$, $fX$ is
completely generated.  In this case, $X$ lies in some subalgebra of
$\fg$ that is isomorphic to $\slc$ via $\phi_\alpha$.

Note that if, with respect to $H_\alpha$, each $x_i$ has even
weight or $f$ has even weight then there is nothing to prove; the
result follows from theorem \ref{even}.  Thus we may assume
without loss of generality that $x_1$ has odd weight $\lambda$ and
that $f = x_1 \cdot g$ where $g \in \faff (M)$ has even weight
$2k$. Moreover, we note that it suffices to assume that
$X=E_\alpha,\ H_\alpha$ or $F_\alpha$.

\noindent Case 1:  ($X=H_\alpha$)  Then
$$[gH_\alpha,x_1H_\alpha] = (\lambda - 2k)x_1 g H_\alpha = (\lambda -
2k) f H_\alpha.$$ Since $gH_\alpha$ has even weight, it is generated
by shears. The result follows from the fact that $\lambda - 2k$ is
odd and hence not zero.

\noindent
Case 2:  ($X=E_\alpha$)  Then
\begin{eqnarray*}
[gH_\alpha,x_1E_\alpha]&=&gH_\alpha(x_1)E_\alpha+f[H_\alpha,E_\alpha]
-x_1E_\alpha(g)H_\alpha\\
&=&(\lambda_1+2)fE_\alpha-x_1E_\alpha(g)H_\alpha.
\end{eqnarray*}
As in Case 1, $gH_\alpha$ is generated by shears.
By Case 1, $x_1E_\alpha(g)H_\alpha$ is
completely generated.  this case now follows from the fact that
$\lambda +2$ is odd and hence not zero. The case where $X=F_\alpha$ is
handled in a fashion similar to case 2.\qed

\medskip

Let $M$ be a Stein manifold, which we think of as already embedded
in some $\cn$.  Suppose we are given a Lie algebra $\fg < \holvec
(\cn )$ of holomorphic vector fields on $\cn$ that are all tangent
to $M$ and that span the tangent space to $M$ at every point.  We
denote by $\faff (M)$ the restriction to $M$ of polynomials in
$\cn$, and write $\vaff (M) := \faff (M) \cdot \fg$, i.e., $\vaff
(M)$ consists of those vector fields that are multiples of the
vector fields in $\fg$ by functions in $\faff(M)$. The next
result, established in \cite{tv}, is a useful sufficiency
criterion for the density property.

\begin{thm}\label{dense}
Suppose the Lie algebra $\fg$ has the following properties.
\begin{enumerate}
\item The Lie algebra $\fg$ consists of linear vector fields (i.e. if
$f$ is a linear function and $\xi \in \fg$, then $\xi f$ is linear)
\item The Lie algebra $\fg$ is semisimple and consists of complete
vector fields,
\item We have $\co (M) \cdot \fg = \holvec (M)$.
\item If $f$ is a linear function and $\xi \in \fg$, then
$f\xi$ is in the Lie algebra generated by complete
vector fields
\end{enumerate}

\noi Then $\vaff (M)$ is generated by complete vector fields.
Thus, in particular, $M$ has the density property.
\end{thm}

\pf Cartan's Theorem A implies that $\vaff (M)$ is a dense subset
of $\holvec (M)$ in the locally uniform topology.  By Theorem
\ref{algebra}, $\C[V] \tensor \fg$ is generated by
$\fg$-overshears, and this remains true in $\vaff (M)$ because
restriction to $M$ preserves Lie brackets. Since the elements of
$\fg$ are complete, so are their overshears, by Proposition
\ref{dror}. \qed

\bs

\noi \rmk We don't ask for $M$ to be parallelized by a basis of
$\fg$; it may happen that $\dim \fg > \dim _{\C} M$.

\bs

As an immediate corollary of Theorem~\ref{dense}, we have

\begin{thm}\label{orbit} Let $G \to {\rm End} (V)$ be a representation
such that $V^* \otimes \fg$ is generated by its $\fg$-overshears.
Then any closed orbit of $G$ has the density property.
\end{thm}

\noi We make use of this theorem by realizing our manifold as a
closed orbit of $G$ via Weyl's Theorem \ref{realization}. (Recall
that a closed submanifold of a Stein manifold is Stein.)

\subsection*{Quadratic vector fields}

In \cite{tv} Theorem \ref{even} was sufficient for establishing
criterion 4 of Theorem \ref{dense}, and thus the density property
on complex semisimple Lie groups. We now proceed to analyze the
situation in greater depth. Recall that criterion 4 requires
vector fields of the form
$$
f\xi, \ \ f \in V^*, \xi \in \fg
$$
to be completely generated. After choosing a basis of $V$ consisting of
weight vectors, and the base $\{E_\alpha,F_\alpha: \alpha \in \Phi^+\}\cup
\{H_\alpha: \alpha \in \Delta\}$ of $\fg$ we are lead to checking
finitely many vector fields.  However, deciding whether any of these vector
fields is completely generated could be fairly difficult. In this section
we investigate to what extent we can reduce the number of vector fields
to be checked so as to obtain the desired result.

To present our sharpening of Theorem~\ref{dense} we introduce the
following.
\begin{defn}
Let $\fg$ be a semisimple Lie algebra. A weight is called extremal if
its orbit under the Weyl group contains a multiple of a fundamental
weight. Given a representation $V$ of $\fg$ a weight vector $v\in V$
is called extremal if the corresponding weight is extremal.
\end{defn}
The term {\it extremal} comes from the fact that the dominant
extremal weights lie in the extremal lines of the dominant
Weyl-chamber (the cone of linear functionals that are positive
with respect to the ordering).

Let $M$ and $\fg$ be as in Theorem \ref{dense}, and fix a Cartan
subalgebra $\fh$ of $\fg$.

\begin{thm}\label{sharpened}
Suppose the Lie algebra $\fg$ has the following properties.
\begin{enumerate}
\item The Lie algebra $\fg$ is semisimple and consists of complete
vector fields that are linear (in the sense of 1 in Theorem
\ref{dense}).

\item We have $\co (M) \cdot \fg = \holvec (M)$.

\item If $f$ is a linear function that is an extremal weight
vector, with extremal weight $\lambda$, and $H_i\in \fh$, is a
simple co-root so that $\lambda(H_i)\neq0$, then $fH_i$ is in the
Lie algebra generated by complete vector fields
\end{enumerate}

\noi Then $\vaff (M)$ is generated by complete vector fields.
Thus, in particular, $M$ has the density property.
\end{thm}

\begin{rmk}
The reduction to $\fh$ is natural from the geometric point of
view: the problematic vector fields are those whose orbit is
$\C^*$. The extremality condition exploits the extra structure of
semi-simple Lie algebras, and does not seem to have a geometric
interpretation.
\end{rmk}

\bigskip
The proof of Theorem \ref{sharpened} is a series of somewhat
technical observations.  We will state these results for general
$\fg$-modules $W$, since they do not rely on any specific
properties of the $\fg$-space at hand, but the reader should keep
in mind that in our application $W$ is the dual of this
$\fg$-space. We begin with the following lemma.

\begin{lem}\label{l-1}
Let $\fg$ be a semisimple Lie algebra, $\fh < \fg$ a Cartan
subalgebra, and $\rho: \fg \rightarrow End(W)$ a representation of
$\fg$ in $W$. Then
$$W \otimes \fg  $$
is generated by $W \otimes \fh$ as $\fg$-module.
\end{lem}

\pf Let $W_0$ be the submodule of $W \otimes \fg$ generated by $W
\otimes \fh$. Given a root $\alpha$ and the corresponding co-root
$H_\alpha$, choose a nonzero  $X_\alpha\in \fg$ so that
$[H_\alpha,X_\alpha]=2X_\alpha$.  The set of all such
$X_{\alpha}$, together with $\fh$, span $\fg$.  If now $w \in W$,
then
$$2 w\otimes X_\alpha = (X_\alpha w)\otimes H_\alpha - X_\alpha(w
\otimes H_\alpha) \in W_0,$$ as desired. \qed

\medskip

\noi The next lemma shows that the weights for only a few of the
$\slc$-subalgebras $\left <  X_{\alpha}, H_{\alpha} \right >$ need
to be even.

\begin{lem} \label{l-2} Suppose $\rho: \fg \rightarrow End(W)$
is a representation of a semisimple Lie algebra $\fg$, and $\fh$
is a Cartan subalgebra of $\fg$. Let $\lambda$ be a weight of
$\fh$ in $W$ and $W_\lambda = \{ w \in W\ ;\ Hw = \lambda(H)w\}$
the corresponding subspace of $W$.  If there exists a co-root
$H_\alpha$ so that $\lambda (H_ \alpha)$ is even and nonzero, then
$W_\lambda \tensor \fh$ is contained in the submodule generated by
the $\fg$-overshears.
\end{lem}

\pf Let $\alpha$ be so that $\lambda(H_\alpha)$ is positive and
even.  Since $\lambda (H_{\alpha}) \neq 0$, we have that
$$\fh= \C H_\alpha \oplus \ker \lambda.$$
Let $W'$ be the smallest subspace of $W$ that is invariant under
$\phi_\alpha(\slc)$ and contains $W_\lambda$. Under the action of
$\phi_\alpha(\slc)$, $\,W'$ has only even weights, and thus by
Proposition \ref{even} it is generated by the $\phi_\alpha(\slc)$
overshears.  These are automatically $\fg$ overshears, showing
that
$$W_\lambda \otimes \C H_\alpha$$
is generated by $\fg$ overshears.  On the other hand, the elements
of
$$W_\lambda \otimes \ker \lambda$$
are overshears for the trivial reason that, for $H \in \ker
\lambda$ and $w \in W_\lambda$, $Hw=0$.
\qed

\begin{lem}\label{l-3} Let $\lambda$ be a weight that is not
a multiple of a fundamental weight. Then there exists a co-root
$H_\alpha$ such that $\lambda(H_\alpha)$ is positive and even.
\end{lem}

\pf Without loss of generality we may assume that $\lambda$ is
dominant. We will use the notation of the Planches (Tables) of
\cite{b}.  We let $\Phi$ denote the root system of $\fg$. For a root
$\alpha$, let $H_\alpha$ be the corresponding co-root in $\fh$.  The
set $\{ H_\alpha \ ;\ \alpha \in \Phi \}$ is a root system dual to
$\Phi$. For a simple root $\alpha_i$, we denote $H_{\alpha_i}$ by
$H_i$. With $\lambda$ as in the hypothesis, one may choose simple
roots $\alpha_i, \alpha_j$, so that $\lambda (H_i)>0,
\lambda(H_j)>0$. Choosing $i,j$ minimal with respect to this property,
we have that $\lambda(H_k)=0$, for all $ i<k<j$.  Assume that
$\lambda(H_i)$ and $\lambda(H_j)$ are odd.  No matter which root
system is being considered, $H_{ij}=\sum_{i\leq k \leq j} H_k$ is a
co-root (although it might not be the co-root corresponding to
$\sum_{i\leq k \leq j}\alpha_k$). Then
$$\lambda(H_{ij})=\lambda(H_i)+\lambda(H_j),$$
so that $\lambda (H_{ij})$ must be positive and even. \qed

\bigskip

\noindent {\it Proof of Theorem~\ref{sharpened}:} By
Lemma~\ref{l-1} it is enough to show that $fH$ is completely
generated for $f \in V^*$, $H\in \fh$. We may assume without any
loss of generality that $f$ is a simultaneous eigenvector of $\fh$
with weight $\lambda$, and that $H=H_\alpha$, for some $\alpha \in
\Delta$. If $\lambda$ is extremal then $fH$ is either complete (if
$\lambda(H)=0$) or else, by assumpion 3 in Theorem
\ref{sharpened}, it is completely generated . We may therefore
assume that $\lambda$ is not a multiple of a fundamental weight.
Lemma~\ref{l-3} shows that there is a co-root $H' \in \fh$ so that
$\lambda(H')$ is even. By Lemma~\ref{l-2} the vector fields in
$V^*_\lambda\otimes \fh $ are completely generated in view of
Proposition~\ref{dror}. Therefore Condition 4 of
Theorem~\ref{sharpened} implies Condition 4 of
Theorem~\ref{dense}, and the result follows.  \qed

\section{Examples}\label{examples-section}

In this section we provide simplified proofs of the density
property for some special examples.  Interestingly, some of these
examples lie beyond the scope of Theorem~\ref{adjoint}.

\subsection*{Small representations}

\begin{defn}
A representation $V$ of $\slc$ is said to be {\emph bounded by 2}
if every weight $n$ satisfies $|n|\le 2$.  More generally, a
representation $V$ of a complex semisimple Lie group $G$ is
bounded by 2 if for any root $\alpha$, the restriction of $V$ to
$\phi _{\alpha} (\slc)$ is bounded by 2.
\end{defn}

\rmk  {\it $(i)$ The notion of boundedness by 2 is independent of the
choice of Cartan subalgebra. $(ii)$  Since the weights of $V^*$ are the
negatives of the weights of $V$, $V^*$ is also bounded by 2.}

\medskip

The next proposition shows that one can use representations
bounded by 2 to produce examples of affine homogeneous spaces
with the density property.  In \cite{tv} we established that if
$V$ is a representation bounded by 2, then $V \otimes \fg$ is
generated by its $\fg$-overshears.  In view of Theorem
\ref{dense}, we therefore have the following proposition.

\begin{prop} \label{std}
If $G$ be a semisimple Lie group admitting a representation $V$
that is bounded by 2, then every closed $G$-orbit in $V$ has the
density property.
\end{prop}

Proposition \ref{std} can be used to establish the density
property in many cases.  We now proceed to do this in several
examples.

\subsection*{Affine quadrics}

Consider the smooth affine subvariety of $\C ^{n+1}$ given by
$$Q_n = \left \{ x_0^2 + ... + x_n ^2 = 1 \right \}.$$

\begin{thm}\label{best-thm}
For any $n \ge 2$, the affine quadric $Q_n$ has the density
property.
\end{thm}

This follows from Proposition \ref{std} and the fact that $Q_n$ is a
closed orbit\footnote{Note, that this homogeneous space is not of adjoint
type.} of the standard representation
$$SO(n+1,\C) \to {\rm End} (\C ^{n+1}),$$
which is bounded by 1.

We will describe this example in some detail. The $n(n+1)/2$ vector
fields $X_{ij},$ $0 \le i < j \le n$, given by
$$X_{ij} := x_i \di _j - x_j \di _i, \qquad 0 \le i < j \le n,$$
are tangent to $Q_n$ and generate the module $\vaff (Q_n)$.  The
relation $$\left [ X_{ij}, X_{kl}\right ] = \delta _{jk} X_{il} -
\delta _{ik} X_{jl} + \delta _{il} X_{jk} - \delta_{jl} X_{ik},$$
where $\delta _{ij}$ is the Kronecker delta, shows that the Lie
algebra
$$\fg = {\rm span} \ \{ X_{ij}\ ;\ 0 \le i < j \le n \}$$
is isomorphic to $\fs \fo (n+1,\C)$.  Let $k$ be the largest
integer such that $2k<n$.  A Cartan subalgebra of $\fg$ is
$${\rm span} \left \{ \ii X_{0,1},..., \ii X_{2k,2k+1} \
\right \}$$ and the weight vectors with nonzero weights are
$$x_0 \pm \ii x_{1},...,x_{2k} \pm \ii x_{2k+1}$$
showing that this representation is bounded by 1.

\subsection*{The space of Lagrangian splittings}

Let $V=\C^{2n}$ with the standard symplectic form
$$\left < x,y \right > := x^tJy,$$
where $x^t$ is the transpose of $x$ and
$$J=\left [
\begin{array}{rr} 0 & -I \\ I & 0 \end{array}
\right ],$$ with the entries $0$ and $I$ signifying the $n \times
n$ zero and identity matrices.  The symplectic group
$$Sp(2n,\C) := \left \{g \in GL(2n,\C )\; :\ g^t J g = J \right \}$$
acts transitively on the Lagrangian subspaces of $V$, i.e., the
subspaces of $V$ that are isotropic for $\left <\ ,\ \right >$
and have maximal dimension.  The set of all Lagrangian subspaces
is a projective manifold, called the Lagrangian Grassmannian, and
thus the stabilizer of any single Lagrangian is parabolic. Choose,
for example, $V_1$ to be the span of the first $n$ standard basis
vectors $e_1=[1,0,0,...,0,0]^t, e_2=[0,1,0,...,0,0]^t, ...,
e_n=[0,0,0,...,0,1]^t$.  The stabilizer of $V_1$ is the subgroup
$$P = \left \{
\begin{bmatrix} a & b \\ 0 & d \end{bmatrix}
\in Sp(2n,\C) \right \} .$$ Consider the set $\cl (V)$ of
Lagrangian decompositions of $V$, i.e., ordered pairs of
Lagrangian subspaces $(L_1,L_2)$, so that $V=L_1 \oplus L_2$. This
set can also be given a holomorphic, and even algebraic structure,
and the resulting manifold $M$ is Stein, supporting a transitive
action by $Sp(2n,\C)$.  If $V_2$ is the subspace generated
$e_{n+1}, ..., e_{2n}$, the stabilizer $$ L=\left \{
\begin{bmatrix} a & 0 \\ 0 & d \end{bmatrix} \in Sp(2n,\C) \right
\} = \left \{ \begin{bmatrix} a & 0 \\ 0 & (a^t)^{-1}
\end{bmatrix}\ ;\ a \in GL(n,\C) \right \}$$ of the pair
$(V_1,V_2)$ is a Levi component of $P$.  Indeed, $L$ is reductive,
$$ U = \left \{ \begin{bmatrix} I & b
\\ 0 & I \end{bmatrix} \in Sp(2n,\C) \right \} $$ is
unipotent, and $P = LU$.

Now $L$ is the centralizer of
$$\lambda= \left [ \begin{array}{rr} I & 0 \\ 0 & -I
\end{array}  \right ] \in M_{2n}(\C),$$
and so $M$ can be identified with the orbit of $\lambda$ in
$M_{2n}(\C)$, the space of $2n \times 2n$ matrices, on which
$Sp(2n,\C)$ acts by conjugation. One checks that $\lambda \in \fs
\fp (2n,\C)$, the Lie algebra of $Sp(2n,\C)$, and so the orbit of
$\lambda$ lies entirely in $\fs \fp (2n,\C)$. The orbit is closed and
is described as
$$
\left\{ \begin{bmatrix} x & y \\ z & -x^t\end{bmatrix}\; :\
y^t=y, \; z^t=z, \; x^2+yz=I \right\}.
$$

Since the adjoint representation of $\fs \fp (2n,\C)$ is bounded by 2,
Proposition \ref{std} implies the following theorem.
\begin{thm}
For every symplectic vector space $V$, the manifold $\cl (V)$ of
all Lagrangian splittings of $V$ has the density property.
\end{thm}

The space of Lagrangian splittings is an example of a {\it
symplectic homogeneous space}, which we now discuss more
generally.

\subsection*{Symplectic homogeneous spaces}

For more details on this brief overview, we refer the reader to
\cite{b}. Let $G$ be a semisimple Lie group. A closed subgroup $P$
is called parabolic if $G/P$ is projective. We give a very brief
description of $P$ in terms of the root system of the Lie algebra
of $G$.

Let $\Phi^+$ denote the set of positive roots, and $\Delta$ the
set of simple roots.  Let $\fn$ be the subalgebra generated by
$X_\alpha$, $\alpha \in \Phi^+$, and $\fb$ the subalgebra
generated by $\fn$ and $\fh$.

Let $B$ denote the connected closed subgroup of $G$ associated to
the Lie algebra $\fb$. The quotient $G/B$ is a projective
manifold, and by a famous theorem of Borel, $B$ is a minimal
parabolic subgroup in the sense that every parabolic subgroup is
conjugate to one that contains $B$. Therefore, up to inner
automorphisms, it is enough to describe those subgroups $P$ whose
Lie algebra $\fp$ contains $\fb$. This is done as follows. Since
$\fp$ contains $\fh$, there is basis of $\fp$ consisting of
eigenvectors of $\fh$, and the set of negative roots admitted by
$\fp$, i.e., the set
$$\left \{ \alpha \in \Phi\ ;\ \exists X_{\alpha} \in \fp \ {\rm
such\ that}\ \forall H\in \fh ,\ [H,X_{\alpha}]= \alpha (H)
X_{\alpha} \ {\rm and}\ - \alpha \in \Phi ^+ \right \},$$ must be
closed under addition. Let $H$ be an element of $\fh$ such that
$\alpha(H)\geq 0$ for all $\alpha \in \Phi^+$, and let $\Phi'=\{
\alpha \in \Phi : \alpha(H)=0 \}$. Then the subspace with basis
$$\left \{X_\alpha : \alpha \in \Phi' \right \} \cup \left \{
H_\alpha : \alpha \in \Delta \right \} \cup \left \{ X_\alpha :
\alpha \in \Phi^+ \backslash \Phi' \right \}$$ is a subalgebra,
and so it is the Lie algebra of a parabolic subgroup. This is the
only way parabolic subalgebras can arise. Moreover, when $H\neq
0$, the algebra constructed here is proper. In that case $P=LU$,
where the Lie algebra of $L$ is generated by
$$\left \{ X_\alpha : \alpha \in \Phi' \right \} \cup \left \{
H_\alpha : \alpha \in \Delta \right \},$$ and thus $L$ is
reductive.  The Lie algebra of $U$ is generated by $\{X_\alpha :
\alpha \in \Phi^+\backslash \Phi' \}$, showing that $U$ is
unipotent. Therefore the Levi component $L$ is the centralizer of
the semisimple element $H$ used in the construction of $\fp$.

For the rest of this section we view $G$ as a subgroup of $End
(\fg)$ through the adjoint representation.  Let $P=LU$ be a
parabolic subgroup, where $L$ arises as the stabilizer of a
semisimple element $H \in \fh$. Since, at the Lie algebra level,
our representation is the derivative of the group action at the
identity, the Lie subgroup of $G$ with Lie algebra equal to the
centralizer of $H$ is the stabilizer of $H$ in the adjoint action.
Therefore $M:=G/L$ is realized as the orbit of $H$ in $\fg$ under
$Ad$, and it is well-known that the adjoint orbit of a semisimple
element is an algebraic submanifold of $\fg$ \cite{b}.

\begin{thm}\label{borel-construction}
The closed adjoint orbit $M=G/L \subset \fg$ described above has
the density property.
\end{thm}

To prove Theorem \ref{borel-construction} it is enough to show
that the conditions of Theorem \ref{dense} are satisfied.  Note
that if $\fg$ has no factors of type $G_2$, the adjoint
representation is bounded by 2 and we could apply Proposition
\ref{std} to these cases. However, the following result covers all
the cases and employs a uniform argument that avoids using the
classification of complex semisimple Lie algebras.  The techniques
are modifications of those used in the proof of Theorem 6.3 in
\cite{tv}.

\begin{thm}\label{adjoints}
The adjoint representation is $\fg$-completely generated.
\end{thm}

\pf Let $G$ be the adjoint group with Lie algebra $\fg$.  In view
of Theorem \ref{algebra}, it suffices to show that for each $x \in
\fg$ and $\vp \in \fg ^*$, $\vp \tensor x$ is generated by
overshears.

Recall that for a representation $V$, the action of $\fg$ on the
dual $V^*$ is given by $x\vp (v) = -\vp (xv)$.

Because $\fg$ is semisimple, it admits a nondegenerate Killing
form, and thus an isomorphism between the adjoint representation
and its dual.  Using the notation of section
\ref{background-section}, this isomorphism is given explicitly as
follows.  Let
$$e_{\alpha} (x) = B(E_{\alpha} ,x)\ \text {for all } x \in \fg.$$
and define $f_{\alpha}$ and $g_{\alpha}$ similarly.

Without loss of generality, we may assume that $\vp $ is any of
the $e_{\alpha},\ h_{\alpha}, f_{\alpha}$ where $\alpha \in \Phi
^+$. Moreover, because of the symmetry between $e_{\alpha}$ and
$f_{\alpha}$, it suffices to prove only that $h_{\alpha} \tensor
E_{\beta},\ h_{\alpha} \tensor H_{\beta},\ f_{\alpha} \tensor
E_{\beta},\ e_{\alpha} \tensor H_{\beta}\ {\rm and}\ e_{\alpha}
\tensor E_{\beta}$ are generated by overshears.  In what follows,
$n$ will be used to denote some integer that may vary from case
to case.

\begin{enumerate}
\item[1:]
$h_{\alpha} \tensor H_{\beta}$ is a shear and $h_{\alpha} \tensor
E_{\beta}$ an overshear.

\noi \pf For the first case,
$$H_{\beta}h_{\alpha} (x) = -h_{\beta} ([H_{\alpha},x]) = -B(H_{\beta}
,[H_{\alpha},x]) = B([H_{\alpha},H_{\beta}],x) = 0.$$ For the
second case,
$$E_{\beta} ^2 h_{\alpha} (x) = - B(H_{\alpha},[E_{\beta},
[E_{\beta},x]]) = - B([[H_{\alpha},E_{\beta}], E_{\beta}],x)=0.$$

\item[2:]
$f_{\alpha} \tensor E_{\beta} = \frac{1}{2}\left( F_{\alpha} (
h_{\alpha} \tensor E_{\beta} ) + n h_{\alpha} \tensor H_{\beta}
\right) ,$ so in view of case 1, $f_{\alpha} \tensor E_{\beta}$ is
generated by overshears.

\item[3:]
$e_{\alpha} \tensor H_{\beta} = -\frac{1}{2}\left ( E_{\alpha} (
h_{\alpha} \tensor H_{\beta} ) + n h_{\alpha} \tensor E_{\beta}
\right ) ,$ so in view of case 1, $e_{\alpha} \tensor H_{\beta}$
is generated by overshears.

\item[4:]
$e_{\alpha} \tensor E_{\beta} = -\frac{1}{2}\left (
E_{\beta}(e_{\alpha} \tensor H_{\beta} ) + c e_{\alpha + \beta}
\tensor H_{\beta} \right )$ for some constant $c$, so in view of
cases 1 and 3, $e_{\alpha} \tensor E_{\beta}$ is generated by
overshears.
\end{enumerate}
This completes the proof.\qed

\medskip

It is known \cite{gs} that if $X$ is a semisimple homogeneous
space that admits a symplectic form, then there is a semisimple
homogeneous space of the type $G/L$ described above, and a finite
holomorphic covering map $\pi : X \to G/L$.  Thus one can apply
Theorem \ref{lift} to obtain the following corollary, which we
believe is worth stating separately.

\begin{cor}\label{symp}
Let $X$ be a semisimple homogeneous space that admits a symplectic
form  and $Y$ a complex manifold with ${\rm dim}(Y) < {\rm dim}
(X)$ such that there exists a proper holomorphic embedding $j: Y
\emb X$. Then there exists another proper holomorphic embedding
$j': Y \emb X$ such that for any $\vp \in {\rm Diff} _{\co} (X)$,
$\vp \circ j (Y) \neq j' (Y)$.
\end{cor}

\section{Proof of Theorem \ref{adjoint}.}
\label{adjoint-proof-section}

In this section we establish Theorem \ref{adjoint} using Theorem
\ref{sharpened}.  In order to reduce the former to the latter, we
provide proofs for a number of simple but perhaps more esoteric
facts about root systems.

We start with some obvious reductions, whose proofs are left to the
reader.

\begin{lem} \label{5.1}\
\begin{enumerate}
\item[a)] If $V$ and $W$ are $\fg$-spaces such that $V
\otimes \fg$ and $W\otimes \fg $ are generated by their
$\fg$-overshears, then $(V \oplus W) \otimes \fg  $ is generated
by its $\fg$-overshears.

\item[b)] Suppose Lie algebras $\fg_1$ and $\fg _2$, and
$\fg_i$-spaces $V_i$, $i=1,2$, are given, so that $V_i \otimes \fg_i $ are
generated by their $\fg _i$-overshears. Then $(V_1 \tensor V_2)
\tensor (\fg_1 \oplus \fg_2)$ is generated by its $\fg_1 \oplus
\fg_2$-overshears.
\end{enumerate}
\end{lem}

\noindent
The representation $\rho$ of $\fg_1 \oplus \fg_2 $ on $V_1 \otimes V_2$
is given by
$$
\rho(\xi _1 \oplus \xi _2) (v_1 \otimes v_2 )= \rho_1(\xi _1)v_1
\otimes v_2 + v_1 \otimes \rho_2(\xi _2 )v_2.
$$

Using Lemma \ref{5.1}, it suffices to consider the case where
$\fg$ is simple and $V$ is an irreducible $\fg$-module.  To this
end, the main result of this section is the following theorem.

\begin{thm}\label{g-mod-conj}
Let $\fg$ be a simple Lie algebra, and let $V$ be an irreducible
$\fg$-space whose highest weight is in the root lattice of $\fg$.
Then $V \otimes \fg $ is generated by its $\fg$-overshears.
\end{thm}

\bs

In view of Theorem~\ref{sharpened} the main task is to handle the
case of multiples of the fundamental weights. Such weight are
sparse in the root lattice, but they do arise.  In the case of
representations with highest weight in the root lattice, the
situation is dealt with by making use of the following facts.
(Again, we use the notation of the Planches (Tables) of \cite{b}.)

\begin{lem}\label{m-omega}
Let $\Phi$ be a simple root system, $B=\{ \alpha_1,
...,\alpha_l\}$ a set of simple roots in $\Phi$, and $\{\omega_1,
..., \omega_l\}$ the corresponding set of fundamental weights. Let
$m\neq 0$ be such that $\lambda=m \omega_i$ is in the root
lattice.
\begin{enumerate}
\item[a)] Assume that $\Phi$ is of type $A_{2l+1}, B_l, C_l, D_l,
E_7, E_8, F_4$, or $G_2$ and that $\omega_i$ is arbitrary.  Then
there exists some root $\alpha$, so that $\lambda(H_\alpha)$ is a
nonzero even integer.

\item[b)] Assume that $\Phi= E_6 $, and that $\lambda=m\omega_i$
for $i=2,3,4,5,6 $. Then there exists some root $\alpha$, so that
$\lambda(H_\alpha)$ is a nonzero even integer.

\item[c)] If $\Phi=A_{2l}$ and $\lambda=m\omega_i$ for some
$i$, or if $\Phi =E_6$ and $\lambda = m \omega_1$, then $m
\omega_i \pm \alpha_i$ is not an integral multiple of a
fundamental weight.
\end{enumerate}
\end{lem}

\pf \underline{a) and b)}:  Recall that for any $\lambda \in
\fh^*$,
$$
\lambda(H_\alpha)= \frac{2 (\lambda,\alpha)}{(\alpha,\alpha)},
$$

\noindent where $(\cdot,\cdot)$ is some multiple of the Killing
form.  To simplify the computations, when working with the
Planches (Tables) of \cite{b}, we will use the ordinary scalar
product in $\R^n$, which in all cases is a scalar multiple of the
Killing form.

\begin{itemize}

\item When $\Phi =A_{2l+1}$, and $\lambda = m \omega_i$ is in
the root lattice, then $(2l+2) | m$, and so $m$ is even.

\item When $\Phi$ is of type $B_l$, then $\omega_1(H_{e_1})=2$
and, for all $i > 1$, $\omega_i(H_{e_1+e_2})=2$.

\item When $\Phi $ is of type $ C_l  $ then, for all $i > 1$,
$\omega_i(H_{e_1+e_2})=2$.  When $i=1$, $\omega_1 \notin \Phi$,
but $2 \omega_1 \in R$, forcing $m$ to be even.

\item When $\Phi$ is of type $ D_l  $ and $ 1 < i < l-1$, $e_1+e_2$
still works for $\omega_i$.  None of $\omega_1, \omega_{l-1},
\omega_l$ is in $\Phi$, and their index in the root lattice is
either 2 or 4, forcing $m$ to be even if $m\omega_1,
m\omega_{l-1}$, or $m \omega_l$ is to be in the root lattice.

\item When $\Phi$ is of type $ E_6, E_7, E_8  $, or $G_2$ we use the
fact that
$$\check{\Phi}  = \left \{ \frac{2\alpha}{(\alpha,\alpha) }: \alpha
\in \Phi \right \} = \Phi.$$ First assume that the pair
$(\Phi,\omega_i)$ is different from the pairs $(E_6, \omega_6)$
and $(E_7, \omega_7)$. Then one finds roots $\alpha = \sum_j c_j
\alpha_j$ so that $c_i$ is even. For $(E_6,\omega_6)$, note that
$m$ must be divisible by 3, and that $\frac{1}{2}(\sum_{j=1}^{5}
e_j - e_6 -e_7 + e_8)$ works. When $\Phi=E_7$, once again $m
\omega_7$ will be in the root lattice only for even $m$.

\item Finally, when $\Phi$ is of type $F_4$, $e_1$ works when
$i=1,2$ and $4$, and $\frac{1}{2}(e_1+e_2+e_3-e_4)$, when $i=3$.
\end{itemize}

\noi \underline{c)}:  This follows from the Euclidean geometry of
the root systems $A_{2\ell}$ and $E_6$. One could also consult the
relevant Planches (Tables) in \cite{b}. \qed

\medskip

The last detail we need is the following generalization of Theorem
\ref{even}.

\begin{prop}\label{old-sl-2}
Let $V$ be an irreducible representation of $\slc$, and suppose $v
\in V$ is such that $Hv=\lambda v$ for some $\lambda \neq 0$. If
neither $E$ nor $F$ annihilates $v$, then $V \tensor \slc$ is
generated by $v \otimes H$ and the $\slc$-overshears of $V \tensor
\slc$.
\end{prop}

\pf If $j$ is the largest integer so that $E^jv \neq 0$, then
$E^jv \otimes E$ and $E^{j-1} v \otimes E$ are overshears, and so
$F^{j+1}(E^{j}v \otimes E)$ and $F^j(E^{j-1}v \otimes E )$ are in
$V$.  Since $FE^jv = (\lambda + 2j) E^{j-1}v$, one can easily
prove, using induction, that
$$ F^{j+1}(E^{j}v \otimes E) = (F^{j+1}E^jv) \otimes E - (j+1) (F^{j}
E^jv )\otimes H +j(j+1) F^{j-1}E^j v \otimes F$$ and
$$ F^{j}(E^{j-1}v \otimes E) = (F^{j}E^{j-1}v) \otimes E -j (F^{j-1}
E^{j-1}v )\otimes H -j(j-1) F^{j-2}E^{j-1}v \otimes F.$$  These
identities make sense, as $j \ge 2$.  It follows that the three
vectors
$$v\tensor H,\ F^{j+1}(E^{j}v \otimes E)\ {\rm and}\
F^{j}(E^{j-1}v \otimes E)$$ are linearly independent. Since $V
\otimes \slc $ is generated by its 3 dimensional $\lambda$
eigenspace, the proof is complete. \qed

\medskip

\noi {\it Proof of Theorem~\ref{g-mod-conj}}.  Let $W$ be the
submodule of $V \otimes \fg$ generated by $\fg$-overshears. By
Lemma~\ref{l-1} we need to show that when $H \in \fh$ and $v$ is a
weight vector of weight $\lambda$,  $v\otimes H \in W$. By
Lemmas~\ref{l-2} and~\ref{l-3} we know that this is true if
$\lambda$ is not an integral multiple of a fundamental weight. The
same arguments show that $ v \otimes H \in W$ even when
$\lambda=m\omega_i$ when $m=0$, or when $\omega_i$ is one of the
fundamental weights listed in Lemma \ref{m-omega}(a).

In the remaining cases, we can assume that $H$ is one of the
co-roots corresponding to a simple root; these were denoted $H_i$
above.  Consider the subalgebra generated by the elements
$E_i,F_i,H_i$, which we denote $\phi_i(\slc)$. Let $M$ be the
smallest subspace of $V$ that contains $V_\lambda$ and is
invariant under $\phi_i(\slc)$, and write $M= \oplus_{\nu} M_{\nu}
$, where each $M_\nu$ is irreducible as a $\phi_i(\slc)$-module.
By complete reducibility, $V_\lambda = \oplus_{\nu} M_\nu \cap
V_\lambda$. Choose $v_\nu \in V_\lambda \cap M_\nu$.

If both $E_i^2v _{\nu}$ and $F_i^2v _{\nu}$ are zero, but $E_i
v_{\nu}$ and $F_i v_{\nu}$ are nonzero, then the representation
$M_\nu$ is isomorphic to the adjoint representation of $\slc$, and
so $H_iv_\nu = 0$, i.e., $m=0$.  When $E_i v_\nu=0$, $v_\nu
\otimes H_i= F_i (E_i v_{\nu} \otimes E_i) - v_{\nu} \otimes E_i
\in W$, and a similar computation works if $F_iv_{\nu} =0$.
Therefore, we may assume that either $E_i^2v_\nu \not = 0$ or $F
_i ^2 v_\nu \neq 0$.  Moreover, by symmetry, it suffices to assume
that $E_i^2 v_\nu \neq 0$.  In this case, since $\lambda+\alpha_i$
is not a multiple of a fundamental weight (Lemma \ref{m-omega}),
we can apply Proposition \ref{old-sl-2} with $v=E_i v_{\nu}$,  and
conclude from Lemmas \ref{l-2} and \ref{l-3} that
$V_{\lambda+\alpha_i} \otimes H \subset W$, as desired.\qed

\medskip

\noi Taking into account the remark at the end of Section
\ref{background-section} and the remark following Theorem
\ref{algebra}, Theorem \ref{adjoint} now follows from Theorems
\ref{orbit}, \ref{realization} and \ref{g-mod-conj}.\qed

\section{Lifting arguments$-$ proof of theorem\ref{lift}}
\label{lift-section}

In the proof of Theorem \ref{lift} we make use of two ideas, due
respectively to Winkelmann \cite{w} and to Forstneri\v
c-Globevnik-Rosay \cite{fgr}.  We shall describe the needed
versions of these ideas below, but we will be brief with certain
parts of the proof, as the proofs that appear elsewhere can be
modified to produce the facts that we need here.

\medskip

We say that two discrete sets $S$ and $T$ on a complex manifold
$M$ are equivalent if there exists $f \in {\rm Diff}_{\co} (M)$
such that $f(S) = T$.  If no such diffeomorphism exists, we say
that $S$ and $T$ are inequivalent. A theorem of Winkelmann
\cite{w} asserts the existence of inequivalent sets on any Stein
manifold.

We begin with the following proposition, which can be established
easily from Winkelmann's Theorem and a counting argument.
\begin{prop}\label{non-tame-for-covers}
Suppose that $\pi : \tilde X \to X$ is a covering map between
Stein manifolds. Then there exist discrete sets $S,T \subset
\tilde X$ that are not equivalent, such that $\pi |S$ and $\pi |T$
are one-to-one.
\end{prop}

Following the method of \cite{fgr}, one needs next to prove the
following result.

\begin{thm}\label{prep}
Let $\pi : \tilde X \to X$ be a finite normal covering map of
Stein manifolds with deck group $\Gamma$, and suppose $X$ has the
density property. Suppose $K \subset \tilde X$ is compact and
holomorphically convex, $A \subset K$ is finite, $p,q \in \tilde
X$ are points not in any of the translates $\gamma K$  of $K$,
$(\gamma \in \Gamma)$ and $\epsilon
>0$.  Then there exists $F \in {\rm Diff} _{\co} (\tilde X)$ such that
\begin{enumerate}
\item $F(p) = q$,
\item $F|A = id$, and
\item $\sup _{x \in K} dist(F(x),x)< \epsilon$.
\end{enumerate}
\end{thm}

The proof of Theorem \ref{prep} is based on Theorem \ref{euler}.
We start with the construction of a vector field on a cover from a
vector field on the base, in such a way that completeness is
inherited. In fact, this is the situation with the usual lifting
of maps in covering space theory: $d\pi : T _{\tilde X} \to T _X$
is also a covering space, and the section $\xi: X \to T _X$ lifts
to a map $\xi ': X \to T _{\tilde X}$ because there are no
topological obstructions.  The map $\xi '$ then factors through
$\pi$, and the factor $\tilde \xi : X \to T_X$ is actually a
section. Rather than proving all of this topologically, we shall
construct the vector field $\tilde \xi$ directly in the case where
the cover is normal. This is the only case needed in our setting.

To this end, let $p \in \tilde X$. Suppose that $\pi ^{-1} \pi p =
\{p, q_2,...,q_d\}$, and let $U_1,...,U_d$ be neighborhoods of $p,
q_2,...,q_d$, respectively, that are mutually disjoint and such
that $\pi |U_j : U_j \to U$ is a diffeomorphism. (In particular,
$\pi (U_i) = \pi (U_j)$.)  We set
$$\tilde \xi _p:= d(\pi |U_1)^{-1} \xi _{\pi p}$$
We leave it to the reader to check that this vector field is well
defined.  Let $\Gamma$ be the deck group of the covering $\pi$.
Since for any $\gamma \in \Gamma$ one has $\pi \gamma = \pi$, the
vector field $\tilde \xi$ so defined is invariant, and thus is the
lift mentioned above.

\begin{lem}\label{lifting-cplt-vflds}
Let $\xi$ be a holomorphic vector field on $X$.  If $\xi$ is
sufficiently close to zero on a compact set $K \subset X$, then
$\tilde \xi$ is small on the compact set $\pi ^{-1} K$. Moreover, if
$\xi$ is complete, then the vector field $\tilde \xi$ is complete on
$\tilde X$.
\end{lem}

\pf The size assertion is obvious.  To see the completeness
assertion, note that, in fact, one can lift every integral curve
of $V$ to $M$.  It is then an easy calculation to see that the
lifted curve is an integral curve of $\tilde V$. Thus the latter
is complete. \qed

\bs

\noi {\it Proof of Theorem \ref{prep}}:  We may replace $K$ by the
union of its translates under the action of the deck group. Let
$\eta_t$ be a time dependent vector field defined on $K \cup \{
p\}$, which is zero on $K$, and whose flow, defined up to time
$1+\delta$, carries the point $p$ to the point $q$ through a path
$c : [0,1]\to \tilde X$ that does not meet $K$. Since, by a
theorem of Stolzenberg, $\cup _{t \in [0,1]} K \times \{t\}$ is
holomorphically convex in $\tilde X \times \C$, we can approximate
$\eta_t$ by a vector field $\xi _t$ on all of $\tilde X$, that is
small on $K$ and vanishes on $A$, and whose flow carries $p$
arbitrarily close to $q$.  We can also assume, after averaging
$\xi _t$ by $\Gamma$, that $\gamma _* \xi _t = \xi _t$ for all
$\gamma \in \Gamma$. This means that we can push $\xi _t$ down to
$X$, where it is a curve in the Lie algebra $\holvec (X)$.  Since
the latter is completely generated, we can approximate $\pi _*
\xi_t$ by a sum of iterated Lie brackets of complete vector
fields. If we now lift all of these complete vector fields to the
cover $\tilde X$, then by Lemma \ref{lifting-cplt-vflds} and the
functoriality of $\holvec$ with respect to mappings we see that
$\xi _t$ is completely generated. An application of Theorem
\ref{euler} finishes the proof.\qed

\bs

\noi Theorem \ref{lift} is an immediate consequence of Proposition
\ref{non-tame-for-covers} and the following theorem.

\begin{thm}\label{interpolate}
Let $f: Y \emb \tilde X$ be a proper holomorphic embedding of
Stein manifolds with ${\rm dim}(Y) < {\rm dim}(\tilde X)$, and let
$S \subset \tilde X$ be any discrete subset such that $\pi |S$ is
$1-1$. Then there exists an embedding $f' :Y \emb \tilde X$ such
that $f'(Y) \supset S$.
\end{thm}

\pf  Let
$$L_0 = \emptyset \subset L_1 \subset {\rm interior} (L_2) \subset
L_2 \subset {\rm interior}(L_3) \subset L_3 \subset ... \subset
\tilde X$$ be a nested family of compact sets such that
$$\tilde X=\bigcup _{\ell \ge 1} L_{\ell} \quad {\rm and} \quad \gamma
L_{\ell} = L _{\ell} \ {\rm for\ all}\ \gamma \in \Gamma.$$ We
decompose $S$ into disjoint finite sets $S_1,S_2,...$ defined by
the requirement that
$$S_j \subset L_j - L_{j-1}.$$
Let $f_0:= f$, and suppose we have obtained an embedding $f_{j-1}
: Y \to \tilde X$ such that
$$f_{j-1} (\Sigma) \supset S_1 \cup ... \cup S_{j-1}.$$
We now apply Theorem \ref{prep} repeatedly to get $f_j$.  To this
end, enumerate
$$S_j = \{ s_{1,j},...,s_{N_j,j} \}.$$
Let $K_{1,j}= L_j$, $A_{1,j} = S_1 \cup ... \cup S_{j-1}$, let $p
\in f_{j-1} \sigma - K_{1,j}$ and $q=s_{1,j}$.  An application of
Theorem \ref{prep} gives us $F_{j,1} \in {\rm Diff} _{\co} (\tilde
X)$ with the properties stated there.  Next, let $2 \le \ell \le
N_j$ and suppose we have obtained $F_{j,\ell -1}$.  Let
$K_{\ell,j} = K_{\ell -1,j} \cup s_{\ell -1,j}$ and $A_{\ell ,j} =
A_{\ell -1,j} \cup s_{\ell -1,j}$, let
$$p \in F_{\ell -1,j} \circ ... \circ F_{1,j} \circ f_{j-1} (\Sigma
)-K_{\ell ,j},$$ and let $q = s_{\ell ,j}$.  Then, again, Theorem
\ref{prep} gives us $F_{\ell,j} \in {\rm Diff} _{\co} (\tilde X)$
with the properties stated there.

We now define
$$f_j := F_{N_j,j} \circ ... \circ F_{1,j} \circ f_{j-1}.$$
It takes some additional care, in the construction of the
$F_{k,j}$, to guarantee that $f_j \to f'$ and that $f'$ is proper.
The details can be carried out in exactly the same way as in
\cite{fgr}.\qed

\bigskip

\' Arp\'ad T\'oth

Department of Analysis

E\"otv\"os Lor\'and University

Budapest, P\'azm\'any P\'eter 1/c

Hungary

email:{\tt toth@cs.elte.hu}

\bigskip

Dror Varolin

Department of Mathematics

University of Illinois at Urbana-Champaign

Urbana, IL 61801  USA

email: {\tt dror@math.uiuc.edu}

\end{document}